\documentclass[a4paper]{article}
\usepackage{amssymb}
\usepackage{amsmath}
\usepackage{latexsym}

\newcommand{\der}[2]{\frac{d#1}{d#2}}
\newcommand{\bin}[2]{\left(  \begin{array}{cc}#1\\#2\\ \end{array}  \right)}

\def\R{{\rm I\!R}}
\def\setN{{\rm I\!N}}

\def\eg {{\sl e.g.}}
\def\setC{\mathbb{C}}

\title{A new basis for eigenmodes  on the Sphere}
\author{M. Lachi{\`e}ze-Rey\\
Service d'Astrophysique, C. E. Saclay\\91191 Gif sur Yvette cedex, 
France\\ 
marclr@cea.fr}
\begin{document}
\maketitle
\abstract{The  usual spherical harmonics $Y_{\ell m}$ form a basis of the vector space  ${\cal V} ^{\ell}$ (of dimension   $2\ell+1$) of the eigenfunctions of the Laplacian on the sphere, with eigenvalue  $\lambda _{\ell} = -\ell ~(\ell +1)$. Here we show the existence of a different  basis  $\Phi ^{\ell}_j$ for ${\cal V} ^{\ell}$, where  $\Phi ^{\ell}_j(X) \equiv  (X \cdot N_j)^{\ell}$, the $\ell ^{th}$ power of the scalar product of the current point with a specific  null vector $N_j$. 
We give explicitly  the transformation properties between the two bases. The simplicity of  calculations  in the new basis allows easy manipulations of the harmonic functions. In particular, we express the transformation rules for the new basis, under any isometry of the sphere.

The development  of  the usual harmonics  $Y_{\ell m}$ into thee new basis (and back) allows   to derive new properties for the $Y_{\ell m}$. In particular, this leads   to  a new relation for the $Y_{\ell m}$, which is a finite version of the well known  integral representation formula. It provides also new development formulae  for the Legendre polynomials and for the  special  Legendre functions.} 

\section{Introduction}

For  a given value of $\ell$,  the $2\ell+1$ spherical harmonics $Y_{\ell m}$ provide a  basis for  the vector space ${\cal V} ^{\ell }$ of   the $2\ell +1$ dimensional  irreductible  representation of the group SO(3).  The same space is  also the  eigenspace   of the Laplacian on $S^2$ with eigenvalue $\lambda _{\ell} = -\ell ~(\ell +1)$. Finally, the reunion of all 
 spherical harmonics  provides a basis for the functions on the sphere. 
 In some sense, they play for the sphere a role analogue to that of the Fourier modes for the plane, and thus appear  very useful for many applications in various fields of mathematics, (classical or quantum) physics, cosmology... 
 
This paper presents the  construction of    another natural basis $\Phi ^{\ell}_j$, $j=-\ell ..\ell$, for ${\cal V} ^{\ell}$, which seems to  have been ignored in the literature.  Each $\Phi ^{\ell}_j$
is  defined as [the reduction to the sphere of] an homogeneous harmonic polynomial in $\R ^3$, which takes the very simple form $(X \cdot N_j)^{\ell}$, where the dot  product extends   the  Euclidean [scalar] dot  product of $\R ^3$ to its  complexification $\setC ^3$, and $N_j$ is a null   vector of $\setC ^3$, that we specify below. 
After defining these functions, we show that they form a basis of  ${\cal V} ^{\ell }$, and we give the explicit  transformation formulae between the two bases. 

The $\Phi ^{\ell}_j$'s have different properties than the $Y_{\ell m}$'s, which  make them more convenient  for particular applications. In particular, they are transformed differently by  the group SO(3), although in a very simple way. Their  main interest comes from  their simplicity, and that of the calculations involving them. In this respect,  we will obtain  new  formulae concerning the usual spherical harmonics. In particular, we obtain    a {\sl finite} version of the usual \emph{integral representation} of the spherical harmonics. And we write  new developments of the Legendre  polynomials and of the special Legendre functions.

\section{Harmonic functions}

We parametrize the [unit] sphere with coordinates $\theta =0..\pi,~\phi=0..2\pi$. An eigenmode [of the Laplacian $\Delta$  in $S^2$]  is  a [complex]   function $f(\theta,\phi)$ on the sphere  which  verifies $\Delta f=\lambda f$.
An   eigenvalue is of the form $\lambda _{\ell}=- \ell~(\ell +1)$, $\ell \in \setN ^{+}$ and  the corresponding eigenfunctions
generate the eigen vectorspace~${\cal V}^{\ell}$, of dimension $2 \ell+1$, which 
  realizes the irreductible unitary representation of the group SO(3).   An usual basis of 
${\cal V}^{\ell}$ is provided  by the [complex]Ê \emph{spherical harmonics} $Y_{\ell m},~m=-\ell..\ell$.
This allows the development of any [complex]Ê function $f$ on the sphere  \begin{equation}
\label{ }
f =\sum _{\ell=0}^{\infty}  \sum _{m=-\ell}^m f_{\ell m}~Y_{\ell m}.
\end{equation}

The unit sphere $S^2$ can be isometrically embedded in $\R ^3$, as the surface $X \cdot X=1$, when a point of $\R ^3$ is written \begin{equation}
\label{ }X=(x,y,z),~
x=r~\cos(\theta),~y=r~\sin (\theta)~\cos(\phi),~z=r~\sin (\theta)~\sin(\phi).
\end{equation} Each function $f \in {\cal V}^{\ell}$ may be seen as the restriction of an harmonic polynomial in $\R ^3$, homogeneous of degree $\ell$, with symmetric and traceless coefficients, that we also note $f$ by an abuse of language: $$f(X)=f(x,y,z)=\sum _{m+n+p=\ell}^{ }  f_{mpq}~x^m~y^n~z^p.$$

\subsection{Complex null vectors}

A complex   vector $Z \equiv (Z^1,Z^2,Z^3)$ is an  element of $\setC ^3$. We extend the Euclidean scalar product in $\R ^3$  to the complex (non Hermitian) inner product  $Z \cdot Z'\equiv \sum_i ~(Z^i)~Ê(Z^{\prime i})$, $i=1,2,3$. A \emph{null}  vector $N$ is defined as having   zero    norm $N \cdot N \equiv \sum _i ~N_i N_i=0 $  (in which  case, it may be   considered as a point on the isotropic cone in $\setC ^3$).
It is straightforward, and well known, that     polynomials  of the form $(X \cdot N)^{\ell}$, homogeneous of degree $\ell$, are  harmonic if and only if $N$ is a null vector. This results from 
$$\Delta _0 (X \cdot N)^{\ell} \equiv  \sum _i  \partial _i ~\partial _i~Ê(X \cdot N)^{\ell}=
\ell ~Ê\sum _i (N _i ~N _i)~Ê(X \cdot N)^{\ell -1}=0,$$ where $\Delta _0$ is the Laplacian of $\R ^3$.
Thus, the    restrictions of  such polynomials are in~${\cal V}^{\ell}$. 
Let us define   the family of null  vectors $N({\alpha})$ with coordinates $(1,i ~\sin  \alpha ,~i ~\cos  \alpha)$, where the   angle $\alpha $ spans the unit circle. 

 It has been shown by \cite{fry} (see also \cite{erd}) that, through this family, the $(k+1)^2$ spherical harmonics of degree $\le k$ of $S^{n-1}$ generate the   spherical harmonics of degree $= k$ of $S^{n}$. In our case   ($n=3$)  this  development  takes the form \begin{equation}
[X \cdot N(\alpha)]^{\ell} =\sum _{m =-\ell}^{\ell} ~ÊÊF_{\ell   m }(X)~Êe^{im~\alpha}.
\label{ }\end{equation}The theorem states that the   $F_{\ell   m }(X)$ are homogeneous harmonic polynomia of degree $\ell$ in $\R ^n$. Multiplying the   previous relation by $e^{-iM~\alpha}$ and   integrating over $\alpha$, we obtain:  \begin{equation} ÊÊF_{\ell   M }(X)=\label{ }\end{equation}
$$\int \frac{d\alpha}{2\pi} ~Êe^{-iM~\alpha}~[X \cdot N(\alpha)]^{\ell} =\int \frac{d\alpha}{2\pi} ~Êe^{-iM~\alpha}~[\cos \theta + i~\sin \theta~\sin (\phi + \alpha )]^{\ell}.$$We recognize the well known integral formula for the spherical harmonics (see, \eg, \cite{niki} p.92), which proves that \begin{equation}
\label{ }
F_{\ell   M }= \frac{1}{2\pi ~B^{\ell}} ~ÊY_{\ell m}, ~~ B^{\ell}_{m} \equiv \frac{1}{4\pi~\ell!}~\sqrt{\frac{(m+\ell)!~(\ell-m)!~(2 \ell+1)}{\pi}}.
\end{equation}
 Thus, we have the development
  \begin{equation}
[X \cdot N(\alpha)]^{\ell} =\frac{1}{2\pi}~\sum _{m =-\ell}^{\ell} ~Ê\frac{Y_{\ell m}(X)}{B^{\ell}_{m}} ~Êe^{im~\alpha},
\label{cinq}\end{equation}which   involves the usual (here, complex) spherical harmonics $Y_{\ell m}$.
Note that the integral formula allows to consider the functions $[X \cdot N_{\alpha}]^{\ell}$, $\alpha \in S^1$ as \emph{coherent states} for the Hilbert space of the harmonic functions (I thank   A. Aribe for this remark). In this sense, all these functions form an overcomplete basis for ${\cal V} ^{\ell }$. We will extract from it a finite basis distinct   from that of the $Y_{\ell m}$. (Note that other analog bases could also be extracted, although less convenient for calculations.)

\subsection{A new basis of eigenmodes on the spheres}

To find a basis of ${\cal V}^{\ell} $ in the form of such polynomials, we  consider the complex $(2 \ell +1)^{th}$ roots  of unity, $\rho ^j$,  $j=-\ell..\ell$,   with $\rho \equiv e^{ i a_{\ell}  }$. 
Their arguments are the  $j~ a_{\ell}  \equiv j~ \frac{2 \pi}{2 \ell+1}$. We will use intensively     the property of the roots of unity, \begin{equation}\label{Dirac}
\sum _j  (\rho^{k})^j=(2 \ell+1)~\delta ^{Dirac} _k,
\end{equation} the Dirac symbol. 
To these roots, we associate (in a given frame)   the $2 \ell +1$ particular null vectors  $N_j \equiv N(j~ a_{\ell} )$, 
~$j=-\ell..\ell$.  
The scalar products are given by $N_j \cdot N_k=1-\cos(\frac{2 \pi ~(j-k)}{2 \ell+1})$.
We want to prove that the functions $\Phi^{\ell}_j:~\Phi^{\ell}_j (X) \equiv (X \cdot N_j )^{\ell}$ form a basis of ${\cal V}^{\ell}$ (hereafter we use the same notation for a function in $\R^3$ and its reduction on the sphere).

Rewriting    equ.(\ref{cinq})  for the     particular value $\alpha=j a$ gives  \begin{equation} \label{phiY} \Phi^{\ell}_j =\frac{1}{2\pi}~\sum _{m =-\ell}^{\ell} ~Ê\frac{Y_{\ell m}}{B^{\ell}_{m}} ~Ê\rho ^{jm}.\end{equation}
 After multiplication  by $\rho^{-jM }$, the  summation over $j$ leads, using (\ref{Dirac}), to an inversion of the formula, namely
\begin{equation}  Y_{\ell m}=\frac{2\pi~B^{\ell}_{m}}{2 \ell +1}~Ê\sum _{j=-\ell}^{\ell} \rho^{-jm}~Ê\Phi^{\ell}_j \label{Yphi}.\end{equation} This  proves that the $\Phi^{\ell}_j$, $\j =-\ell,\ell$ form a basis of ${\cal V}^{\ell} $.
The two last formulae  express the change of bases. 
As an illustration, we   give the first developments in Appendix A.

The basis $(\Phi^{\ell}_j)$ appears not orthogonal. However, it easy to calculate the scalar products\begin{equation}
\label{ }
<\Phi^{\ell}_j, \Phi^{\ell '}_{j'}>:=\frac{\delta _{\ell \ell '}}{(2\ell+1)^2}~\sum _m~\rho^{m (j-j')}~\frac{1}{(D_{\ell m})^2}\end{equation}
$$=\frac{4\pi ~(\ell !)^2}{(2\ell+1)~(2\ell)!} [2~ \cos(\frac{(j-j')~\pi}{2\ell+1})]^{2\ell}$$ 
$$=\frac{4\pi ~(\ell !)^2~2^{\ell}}{~(2\ell+1)~(2\ell)!~2^{2\ell}} [2-N_j \cdot N_{j'}]^{\ell}$$ 

\subsection{New developments of Legendre polynomials and  special Legendre functions}

The    expansion of   (\ref{Yphi}) with the binomial coefficients reads:
$$ Y_{\ell m}(\theta,\phi)=\frac{2\pi~B^{\ell}_{m}}{2 \ell +1}~
 \sum _j ~\rho^{-m~j}~\sum _{p=0}^{\ell}  \bin{\ell}{p}~(\cos \theta)^{\ell-p}~[ \frac{\sin \theta}{2}~(z~\rho^j-\frac{1}{z~\rho^j } )]^{p},$$where  $z \equiv e^{i\phi}$.
We develop again,   permute the summation symbols, and rearrange the terms : the sum takes the form
$$\sum _{p=0}^{\ell}  \bin{\ell}{p}~(\cos \theta)^{\ell-p}~\left( \frac{\sin \theta}{2}\right)^{p}~\sum _{q=0}^{p}  \bin{p}{q} ~(-1)^{p-q} ~z^{2q-p} ~\sum _j ~\rho^{j(2q-p-m)} .$$

Because of     the property (\ref{Dirac}),  the term in the sum  is non zero only when $2q-p-m=0$. 
The term in   the development is non zero only when $p+m$ is even, and $0 \le q \le p \le \ell$.  Taking  these conditions   into account, we   rewrite the sum as 
\begin{equation} \label{zut}  \sum _{q=\max(0,m)}^{[\frac{\ell+m}{2}]} 
\frac{\ell !~2^{m-2q}~   ~(-1)^{q-m}}{(\ell +m-2q)!~Êq!~Ê(q-m)!}
 ~ (\cos \theta)^{\ell+m-2q}~( \sin \theta )^{2q-m},\end{equation}
  where the bracket means entire value.

As expected, $Y_{\ell m}$  is proportional to $z^m$, and  is thus  an eigenfunction 
of the rotation operator $P_x$.
It results  the new development formula for the   associated Legendre functions (defined as  
$P_{\ell}^m(x)Ê\equiv  (1-x ^2)^{m/2} ~\der{P_{\ell}(x)}{x^m}$):\begin{eqnarray}\label{hjk}
P^{\ell} _m (x)=    
   ~\sum  _{q=\max(0,m)}^{[\frac{\ell+m}{2}]} 
\frac{(m+\ell)!~(-1)^{q-m}~2^{m-2q}}{(\ell+m-2q)!~q!~(q-m)!} ~~(1-x^2) ^{q-m/2}~    x^{\ell +m -2q}.\end{eqnarray}

We can specify to  the Legendre polynomials, by putting $m=0$. We obtain
\begin{eqnarray}\label{hjk}
P^{\ell}  (x)=    \ell !~\sum  _{q=0}^{[\frac{\ell }{2}]} 
\frac{ 2^{2q}}{(\ell-2q)!~(q !)^2} ~~(x^2-1) ^{q}~    x^{\ell   -2q}.\end{eqnarray}

\subsection{The integral representation becomes finite}

The  integral representation of the spherical harmonics is   the well known (\cite{niki} p.92):\begin{equation}
\label{ }
Y_{\ell m}(\theta,\phi)=B^{\ell}_{m}~\int _{-\phi}^{2\pi -\phi} d\alpha ~e^{-im \alpha}~[\cos \theta+i~\sin \theta~\sin (\phi +\alpha)]^{\ell},
\end{equation}with $ B^{\ell}_{m} \equiv \frac{1}{4\pi~\ell!}~\sqrt{\frac{(m+\ell)!~(\ell-m)!~(2 \ell+1)}{\pi}}$.
The explicit development of  (\ref{Yphi})    provides a
a {\sl finite version} of the latter, as 
\begin{equation}
\label{ }
Y_{\ell m}(\theta, \phi)=\frac{2\pi~B_{\ell m}}{2 \ell +1} \sum  _{j=-\ell}^{\ell} ~e^{-i m j~\alpha }~ 
 [\cos \theta + i~\sin \theta~\sin (\phi +j~\alpha  )]^{\ell}, ~\alpha  \equiv \frac{2\pi }{2\ell +1}.\end{equation}
 
\subsection{Group action}

The vector space  ${\cal V}^{\ell}$ form an $(2 \ell +1)$ dimensional   IUR of SO(3), $T$,  whose action is defined through \begin{equation}
\label{ }
T:~SO(3) \ni g \mapsto T_g~:~f \mapsto T_g f: T_g f(x) \equiv f(gx).
\end{equation}
This action is completely defined  by the transformation laws for a basis of 
${\cal V}^{\ell}$. For instance, the usual spherical harmonics have the property that they transform very simply (a   multiplication by a complex number) under a selected SO(2) subgroup of SO(3): they have been precisely chosen as eigenfunctions of the angular momentum operator, in addition to the Laplacian (which is the Casimir operator of the group).

The  functions of  the  new basis transform in a different way under the group: developing the transformed function in the basis, we can write     \begin{equation}
\label{ }
T_{g}:~\Phi^{\ell}_j \mapsto T_{g}\Phi^{\ell}_j = \sum _k~ÊG_{j}^{(\ell)k}[g]~\Phi^{\ell}_k, ~\forall g \in \mbox{SO}(3) .
\end{equation}
On the other hand, by definition of the representation,   $ T_{g}\Phi^{\ell}_j (X)=\Phi^{\ell}_j (gX)$. This leads to the relation \begin{equation}
\label{ }
(gX \cdot  N_j)^{\ell}=\sum _k ~ÊG_{j}^{(\ell)k}[g]~(X \cdot N_k)^{\ell}.
\end{equation}  

To go further, we define the two vectors   $ \alpha \equiv (0,1,-i)$ and $\beta \equiv (1,0,0)$, such that $\alpha \cdot N_j=\rho ^j$ and $\beta \cdot N_j=1$, and we apply the formula above to  the vector   $X  = \alpha + t \beta$. Developing and identifying the powers of $t$, we obtain \begin{equation}
\label{ }
(g \alpha \cdot N_j)^{p}~Ê(g \beta \cdot  N_j)^{\ell -p}=\sum _k ~ÊG_{j}^{(\ell)k}[g]~(\alpha \cdot N_k)^{p}~~~~\forall p =0 .. k.
\end{equation}
The  solution  
provides the explicit form of  the coefficients $$ÊG_{j}^{(\ell)k}[g]= \frac{1}{2 \ell +1}~ \sum _{n}    \rho ^{-k  n} ~(g   \alpha \cdot N_j  )^{n}~Ê(g \beta \cdot  N_j)^{\ell -n} ,$$ for an arbitrary $g \in $ SO(3). This provides the transformation properties of any function developed in the basis.

If we select a subgroup SO(2) in SO(3), there is a basis (in which   the pole $[1,0,0]$ remains  fixed) where its elements are given by the matrices 
$g_{\psi}=\begin{pmatrix}
      1&0    &0\\
      0&\cos \psi  &-\sin \psi\\
      0&\sin \psi    &\cos \psi \end{pmatrix},~\psi  \in [0..2\pi].$  
       This implies\\ $g_{\psi}   \alpha \cdot N_j=w \rho ^j$, $g_{\psi}   \beta \cdot N_j=1$, with $w \equiv e^{i \psi}$, and thus     \begin{equation}
 G_{j}^{(\ell)k}=\frac{1}{2 \ell +1}~\sum _{n=-\ell}^{\ell}~w^{n}~\rho^{n(j-k)}.
\end{equation} Given the correspondence (\ref{phiY},\ref{Yphi}),  this  is equivalent to the well known formula
describing the  rotation properties of the spherical harmonics, $Y_{\ell m} \mapsto w^m~ÊY_{\ell m}$.

Note that in the case $w= \rho^K$, a ($2\ell+1$) root of unity, the matrix becomes diagonal:  
$ G_{j}^{(\ell)k}=\frac{1}{2 \ell +1}~\sum _{n=-\ell}^{\ell}~\rho^{nK}~\rho^{n(j-k)}=\delta ^{Dirac}_{K+j-k}$ (a remark due to A. Aribe).
  
\section{Conclusion} 

We have presented a new basis for the eigenfunctions of the Laplacian on the sphere, distinct from the usual basis of spherical harmonics  $Y_{\ell m}$, and with different properties. Its very simple expression  allows easy calculations. In particular, we gave a general formula which gives their transformation properties under the isometry group of the sphere, SO(3). This led, using the transformation formulae, to  easy calculations of the transformation properties  of the usual spherical harmonics, not only for a privilegied group SO(2) of SO(3). In addition,   new formulae were derived for the $Y_{\ell m}$, as well for the Legendre polynomials and the  Legendre special  functions, of great potential use   for calculations with harmonic functions.

Subsequent work  will generalize this construction to $S^3$. Beside the intrinsic interest, this will provide  an explicit formulation of the transformation properties of the harmonic functions on $S^3$ under the isometry group SO(4). In turn, this will allow to select those functions which remain invariant under  specified elements of SO(4). This opens the way to calculate the eigenfunctions [of the Laplacian] on any spherical space $S^3 /\Gamma$, which remain presently unknown in general.

\section{Appendix A} 
We develop the first basis functions in Spherical Harmonics:
 \begin{eqnarray}
\Phi ^{1}_0&=\sqrt{\frac{2\pi}{3}}~&[Y ^{1}_{-1}+\sqrt{2}~Y ^{1}_{0}+Y ^{1}_{1}],\nonumber \\   
\Phi ^{1}_{-1}&= \sqrt{\frac{2\pi}{3}}~&\left[\frac{- 1 +i~ \sqrt3}{2} ~Y ^{1}_{-1} +\sqrt{2}~Y ^{1}_{0}-  \frac{1 +i~ \sqrt3}{2}~Y(1)\right],\nonumber \\ 
\Phi ^{1}_1&= \sqrt{\frac{2\pi}{3}}~&\left[- \frac{1 +i~ \sqrt3}{2}Y ^{1}_{-1} +\sqrt{2}~ Y ^{1}_{0}+  \frac{- 1 +i~ \sqrt3}{2}~Y(1)\right];\nonumber \\ 
\Phi ^{2}_{0}&=
\frac{2\sqrt{\pi}}{\sqrt{5}}~&\left[\frac{Y ^{1}_{-2}}{\sqrt{6}}
+ Y ^{1}_{-1}~\sqrt{\frac{2}{3}}+
 Y ^{1}_{0}
+Y ^{1}_{1}\sqrt{\frac{2}{3}}+\frac{Y ^{1}_{2}}{\sqrt{6}}\right],\nonumber \\ 
\Phi ^{2}_{-1}&=
\frac{2\sqrt{\pi}}{\sqrt{5}}~&\left[\frac{ e^{\frac{4i\pi}{5}}~Y ^{1}_{-2}}{\sqrt{6}}
+ e^{\frac{2i\pi}{5}}~Y ^{1}_{-1}~\sqrt{\frac{2}{3}}+
 Y ^{1}_{0}
+e^{\frac{-2i\pi}{5}}~Y ^{1}_{1}\sqrt{\frac{2}{3}}+\frac{e^{\frac{-4i\pi}{5}}~Y ^{1}_{2}}{\sqrt{6}}\right],\nonumber \\ 
\Phi^{2}_{1}&=
\frac{2\sqrt{\pi}}{\sqrt{5}}~&\left[\frac{ e^{\frac{-4i\pi}{5}}~Y ^{1}_{-2}}{\sqrt{6}}
+ e^{\frac{-2i\pi}{5}}~Y ^{1}_{-1}~\sqrt{\frac{2}{3}}+
 Y ^{1}_{0}
+e^{\frac{2i\pi}{5}}~Y ^{1}_{1}\sqrt{\frac{2}{3}}+\frac{e^{\frac{4i\pi}{5}}~Y ^{1}_{2}}{\sqrt{6}}\right],\nonumber \\ 
 \Phi^{2}_{2}&=
\frac{2\sqrt{\pi}}{\sqrt{5}}~&\left[\frac{ e^{\frac{-8i\pi}{5}}~Y ^{1}_{-2}}{\sqrt{6}}
+ e^{\frac{-4i\pi}{5}}~Y ^{1}_{-1}~\sqrt{\frac{2}{3}}+
 Y ^{1}_{0}
+e^{\frac{4i\pi}{5}}~Y ^{1}_{1}\sqrt{\frac{2}{3}}+\frac{e^{\frac{8i\pi}{5}}~Y ^{1}_{2}}{\sqrt{6}}\right],\nonumber \\ 
 \Phi^{2}_{-2}&=
\frac{2\sqrt{\pi}}{\sqrt{5}}~&\left[\frac{ e^{\frac{8i\pi}{5}}~Y ^{1}_{-2}}{\sqrt{6}}
+ e^{\frac{4i\pi}{5}}~Y ^{1}_{-1}~\sqrt{\frac{2}{3}}+
 Y ^{1}_{0}
+e^{\frac{-4i\pi}{5}}~Y ^{1}_{1}\sqrt{\frac{2}{3}}+\frac{e^{\frac{-8i\pi}{5}}~Y ^{1}_{2}}{\sqrt{6}}\right],\nonumber 
\end{eqnarray}

\end{document}